\documentstyle[12pt]{article}
\newtheorem{theorem}{'heorem}
\newtheorem{lemma}{Lemma}
\newtheorem{note}{Note}

\newtheorem{proposition}{Proposition}
\def\Proof{{\bf Proof}}

\font\Bbb=msbm10 scaled 1200
\def\R{\hbox{\Bbb R}}

\def\C{\hbox{\Bbb C}}
\def\beq{\begin{equation}}
\def\eeq{\end{equation}}

\def\Id{\mathop{\rm Id}}
\def\DV{\mathop{\rm DV}}
\def\ker{\mathop{\rm Ker}}

\begin{document}
""Š 512.542+517.551

\begin{center}
GEOMETRY OF CROSS RATIO \\

\smallskip
Zelikin M.I.
\footnotetext{This research  was supported in part by the Russian Foundation
for Basic Research grant 04-01-00735, and Mintechprom RF grant
NSh-304.2003.1}
\end{center}

\bigskip

Generalization of the cross ratio to polarizations of linear finite
and infinite-dimensional spaces (in particular to Sato Grassmannian)
is given and explored. This cross ratio appears to be a cocycle of the
canonical (tautalogical) bundle over the Grassmannian with
coefficients in the sheaf of its endomorphisms. Operator analog of the
Schwarz differential is defined. Its connections to linear
Hamiltonian systems and Riccati equations are established. These
constructions aim to obtain applications to KP-hierarchy.

\bigskip


\smallskip


\bigskip

         {\bf Sato Grassmanian and KP-hierarchy}

\bigskip

We shall use the following standard notations. As a field of constants
$k$ one may consider both $\R$ or $\C$.
For a basic ring
$B = k[[x]] = \{\sum_{i \ge 0} a_ix^i \,| \, a_i \in k\}$ --
(the formal Taylor series) let
$E = k[[x]]((\partial^{-1}))$ -- be the ring of formal
pseudodifferential operators with coefficients in $B$.

Consider the direct sum decomposition $E = E_+ \oplus E_-$,
where $E_+ = B[\partial]$ is the ring of differential operators,
and $E_- = B[[\partial^{-1}]]\partial^{-1}$ -- is the ring of operators
with the negative order (Volterra type operators). So that for any
$P \in E$ one have $P = P_+ + P_-$, where $P_+$ is a differential part
and $P_-$ is a Volterra part of $P$. A decomposition of an
infinite-dimensional space on a direct sum of its two
infinite-dimensional subspaces shall be called polarization.

Let $H$ be a Hilbert space, equipped with a polarization, i.e.
with a decomposition on two closed orthogonal subspaces
\beq
H = H_+ \oplus H_-.
\label{1}
\eeq

It is convenient to define this decomposition by an operator of
complex structure, i.e. by the unitary operator $J: H \to H$ which is
equal to $+\Id$ on $H_+$ and to $-\Id$ on $H_-$. The General
restricted group $GL_{res}(H)$ is defined as a subgroup of $GL(H)$
consisted of operators $A$ such that the commutator $[J,A]$ is a
Hilbert-Schmidt operator. Equvivalently, if $A \in GL(E)$ is
represented as $(2 \times 2)$-matrix

$$
A =
\left (
\begin{array}{cc}
a & b \\
c & d
\end{array}
\right )
$$
\noindent relative to the decomposition~(\ref{1}), then
$A \in GL_{res}(H)$ iff $b$ and $c$ are Hilbert-Shmidt operators.
It follows that $a$ and $d$ are Fredholm operators.
$U_{res}(H)$ is a subgroup of $GL_{res}(H)$, which consists of unitary
operators.

Sato Grassmannian $SGr(H)$ is the set of all closed subspaces
$W \in H$ such that

1. The orthogonal projection $pr_+: W \to H_+$ is a Fredholm operator

2. The orthogonal projection $pr_-: W \to H_-$ is a
Hilbert-Shmidt operator.

It is known~\cite{PS} that the subgroup $U_{res}(H)$ of the group
$GL_{res}(H)$ acts on $SGr(H)$ transitively, and the stabilizer of the
subspace $H_+$ is $U(H_+) \times U(H_-)$. Hence $SGr(H)$ is represented
as a homogeneous space
$SGr(H) =
U_{res}(H) / ((U(H_+) \times U(H_-))$. $SGr(H)$ is an
infinite-dimensional manifold modelled by Hilbert space
$J_2(H_+,H_-)$ of Hilbert-Shmidt operators $T:H_+ \to H_-$
where the norm $T \in J_2(H_+,H_-)$ is defined as
$(\sum_i \parallel Te_i \parallel^2)^{1/2}$, and $\{e_i\}$ is a full
orthonormal basis on $H_+$. The set of such operators (more exactly, a
chart on $SGr(H)$ corresponding to this set) is called the big cell and
is denoted by $SGr_+$. The group $GL_{res}(H)$ acts on
$SGr(H)$ by M\"obius transforms
$$
\left ( a \quad b \atop c \quad d \right ) \cdot T =
(c + dT)(a + bT)^{-1}.
$$

The tangent space to the manifold $SGr(H)$ at a point $W$ is
$T_W = \mathop{\rm Hom} (W,H/W)$.

Consider the affine subspace $E' = \partial + H_-$
of the linear space $H$ and the multiplicative group
$G = 1 + H_-$. Both of them have the tangent space $H_-$.
Define the mapping $\xi: G \to E'$ for $S \in G$ by the formula
$\xi: S \mapsto L = S \partial S^{-1}$ and the mapping
$\eta: G \to Gr_+$ by the formula $\eta: S \mapsto S^{-1}w_0$,
$w_0 \in SGr_+$ being a fix subspace. Differentials of these mappings
at a point $S$ acts on $A \in H_-$ by formulae
$d \xi : A \mapsto [AS^{-1},L]$ and
$d \eta : A \mapsto -S^{-1}A$ accordingly.

The main result (known as Sato's correspondence)~\cite{PS},
~\cite{S} is: The operator of multiplication by $z^{-n}$
(which induces $S \in \mathop{\rm Hom} (W,H/W)$) transfers by
$d \xi (d \eta)^{-1}$ into the commutator $[L_+^n, L]$, i.e. into the
right hand side of KP equation. Given $n$, the operator $z^{-n}:H \to H$
defines a fix vector field on $SGr(H)$, i.e. the right hand side of
"constant coefficient" Riccati equation~\cite{Z}. Thus KP-hierarchy
corresponds to a countable polysystem of mutually commuting flows of
operator Riccati equations on Sato Grassmannian. Analogues theory
for noncommutative rings of coefficients was developed in~\cite{Mu},
and for multidimensional $x$ in~\cite{P}. To investigate KP-hierarchy
it is supposed to use the following generalization of the cross ratio
~\cite{Ze}.

\bigskip

 {\bf Operator cross-ratio}

\bigskip

Consider first the finite-dimensional case and subspaces of half
dimension. Let ${\cal P}_1,{\cal P}_2,{\cal P}_3,{\cal P}_4$ be four
$n$-dimensional subspaces in $\R^{2n}$ which corresponds to four
points with matrix coordinates $P_1,P_2,P_3,P_4$
of the big cell of the manifold $Gr_n(\R^{2n})$. Suppose that
${\cal P}_1,{\cal P}_2$ and ${\cal P}_3,{\cal P}_4$ define
polarizations of
$\R^{2n}$, i.e. ${\cal P}_1 \oplus {\cal P}_2 =
{\cal P}_3 \oplus {\cal P}_4 = \R^{2n}$. The polarization
${\cal P}_i,{\cal P}_j$ will be denoted by $\Pi_{ij}$.
The class of matrices which are similar to the matrix
$(P_1-P_2)^{-1}(P_2-P_3)(P_3-P_4)^{-1}(P_4-P_1)$
(invers matrices are defined since $\Pi_{12}$ and $\Pi_{34}$
are polarizations) is an invariant of the ordered four points of
Grassmannian relative to M\"obius transformations~\cite{Z}.
It is called matrix cross-ratio.

The projection parallel to a subspace ${\cal P}_i$ will be denoted by
$\pi_i$ or by the figure $i$ above the arrow which gives the
corresponding mapping. If ${\cal P}_i$ and ${\cal P}_j$ defines a
polarization, then the image of $\pi_i$ in ${\cal P}_j$ is uniquely
defined.

\begin{theorem}
\label{theorem1}
Let  $\Pi_{12}$ and $\Pi_{34}$ be polarizations. Then
$\DV({\cal P}_1,{\cal P}_2;{\cal P}_3,{\cal P}_4)$ is the matrix
of the composition mapping
${\cal P}_1 \stackrel{4}{\longrightarrow}{\cal P}_3
 \stackrel{2}{\longrightarrow}{\cal P}_1$
of the space ${\cal P}_1$ on itself.
\end{theorem}

\Proof.
Let $(g,P_3g)$ be the projection of an element $(f,P_1f)$ of the space
${\cal P}_1$ on ${\cal P}_3$ parallel to ${\cal P}_4$. This means that
$(f-g,P_1f-P_3g) \in {\cal P}_4$, i.e. $P_1f-P_3g = P_4(f-g)$.
Hence, $g = (P_4-P_3)^{-1}(P_4-P_1)f$. Similarly, the projection of an
element $(g,P_3g)$ on ${\cal P}_1$ parallel to ${\cal P}_2$ is
$(h,P_1h)$, where $h = (P_2-P_1)^{-1}(P_2-P_3)g =
(P_1-P_2)^{-1}(P_2-P_3)(P_3-P_4)^{-1}(P_4-P_1)f$.
$\Box$

\begin{note}

The $\tau$-function in the theory of integrable systems {\rm
\cite{M}} is in essence only the determinant of an operator
cross-ratio; remaining invariants of the corresponding operators was
not used before. Hence, in {\rm \cite{M}} was proved less general
assertion concerning only about the determinant of the cross ratio.

\end{note}

Our theorem allows to remove the restriction
of half-dimension and to define cross-ratio for a pair of
polarizations of the space $\R^m$ with similar dimensions, i.e.
$\dim {\cal P}_1 = \dim {\cal P}_3 = k; \;
\dim {\cal P}_2 = \dim {\cal P}_4 = m-k.$ Moreover, this gives the
possibility to define an operator cross-ratio for infinite-dimensional
case. In so doing, the invariance of the operator cross-ratio relative
to the M\"obius group is inherited by the construction itself due to
the linearity of projection operators. Note that the composition
mapping is invariant but its matrix is defined up to conjugation

Let ${\cal P}_1, {\cal P}_2, {\cal P}_3, {\cal P}_4$ be subspaces of
linear space $H$. The operator
$\DV({\cal P}_1,{\cal P}_2;{\cal P}_3,{\cal P}_4)$ defined in theorem 1
will be denoted by $\DV_{12;34}$.

We need another expression for the cross ratio. Let
$H=E_1 \oplus E_2$ be a polarization of $H$. We shall call $E_1$ the
horizontal subspace and $E_2$ the vertical one. We suppose that
${\cal P}_1,{\cal P}_3 \in E_1$, ${\cal P}_2,{\cal P}_4 \in E_2$.
Let the subspaces ${\cal P}_1$, ${\cal P}_3$ have as its coordinates
matrices $P_1$ and $P_3$ correspondingly. For the subspaces ${\cal P}_2$
and ${\cal P}_4$ we change the role of the vertical and the horizontal
subspaces. Let in this new system of coordinates the subspaces
${\cal P}_2$ and ${\cal P}_4$ have as its coordinates the matrices
$P_2$ and $P_4$ correspondingly.

\begin{proposition}
\label{proposition1}
Let $\Pi_{12}$ and $\Pi_{34}$ are polarizations and
let the projections parallel to ${\cal P}_2$ and ${\cal P}_4$ are
isomorphisms of the spaces
${\cal P}_1$ and ${\cal P}_3$. Then the matrix
$\DV({\cal P}_1,{\cal P}_2;{\cal P}_3,{\cal P}_4)$
of the composition mapping
${\cal P}_1 \stackrel{4}{\longrightarrow}{\cal P}_3
 \stackrel{2}{\longrightarrow}{\cal P}_1$
of the space ${\cal P}_1$ on itself has the form
$$
(P_2P_1 - I)^{-1}(P_2P_3 - I)(P_4P_3 - I)^{-1}(P_4P_1 - I).
$$
\end{proposition}

\Proof.

Let $(g,P_3g)$ be the projection of an element $(f,P_1f)$ of the space
${\cal P}_1$ on ${\cal P}_3$ parallel to ${\cal P}_4$. This means that
$(f-g,P_1f-P_3g) \in {\cal P}_4$.
As soon as the roles of of the vertical and the horizontal subspaces
interchange for ${\cal P}_4$ we have
$P_4(P_1f-P_3g) = (f-g)$. The matrix $(P_4P_3 - I)$ is invertible
hence we have $g = (P_4P_3 - I)^{-1}(P_4P_1 - I)f$. Analogousely
the projection of the element $(g,P_3g)$ on ${\cal P}_1$ parallel to
${\cal P}_2$ is
$(h,P_1h)$, where $h = (P_2P_1 - I)^{-1}(P_2P_3 - I)g =
(P_2P_1 - I)^{-1}(P_2P_3 - I)(P_4P_3 - I)^{-1}(P_4P_1 - I)f$.

$\Box$

\bigskip

\begin{lemma}
\label{lemma1}
Let projections parallel to ${\cal P}_2$ and ${\cal P}_4$ be
isomorphisms of the spaces ${\cal P}_1$ and ${\cal P}_3$. Then
$\DV_{12;34} = \DV_{34;12}$.
\end{lemma}

\Proof.  Cosider the mappings
$$
\begin{array}{l}
{\DV}_{12;34}:
{\cal P}_1 \stackrel{4}{\longrightarrow}{\cal P}_3
\stackrel{2}{\longrightarrow}{\cal P}_1; \\
{\DV}_{34;12}:
{\cal P}_3 \stackrel{2}{\longrightarrow}{\cal P}_1
\stackrel{4}{\longrightarrow}{\cal P}_3.
\end{array}
$$
If one identifies ${\cal P}_1$ and ${\cal P}_3$ by using projection
$\pi_4$, then both maps will coinside.
$\Box$

Thus the cross-ratio does not depend on the order of polarization
pairs.

\begin{lemma}
\label{lemma2}
Let ${\cal P}_1 \oplus {\cal P}_2 = {\cal P}_3 \oplus {\cal P}_4 =
{\cal P}_1 \oplus {\cal P}_4 = {\cal P}_3 \oplus  {\cal P}_2 = H$.
Then $\DV_{12;34} = \DV_{14;32}^{-1}$.
\end{lemma}

\Proof.  Consider the mappings
$$
\begin{array}{l}
{\DV}_{12;34}:
{\cal P}_1 \stackrel{4}{\longrightarrow}{\cal P}_3
\stackrel{2}{\longrightarrow}{\cal P}_1; \\
{\DV}_{14;32}:
{\cal P}_1 \stackrel{2}{\longrightarrow}{\cal P}_3
\stackrel{4}{\longrightarrow}{\cal P}_1.
\end{array}
$$
We have
${\cal P}_1 \stackrel{4}{\longrightarrow}{\cal P}_3 =
({\cal P}_3 \stackrel{4}{\longrightarrow}{\cal P}_1)^{-1}$ ¨
${\cal P}_1 \stackrel{2}{\longrightarrow}{\cal P}_3 =
({\cal P}_3 \stackrel{2}{\longrightarrow}{\cal P}_1)^{-1}$,
which follows the lemma.
$\Box$

\begin{lemma}
\label{lemma3}
Let ${\cal P}_1 \oplus {\cal P}_2 = {\cal P}_3 \oplus {\cal P}_4 =
H$.
Then $\DV_{12;43} = \Id - \DV_{12;34}$.
\end{lemma}

\Proof. Design the image of $h_1 \in {\cal P}_1$ under the mapping
$\DV_{12;43}: {\cal P}_1 \stackrel{3}{\longrightarrow}{\cal P}_4
\stackrel{2}{\longrightarrow}{\cal P}_1$.
We have $\pi_3 h_1 = h_1 + h_3$, where
$h_3 \in {\cal P}_3$, and $h_1 + h_3 \in {\cal P}_4$. Further on,
$\pi_2 (h_1 + h_3) = h_1 + h_3 + h_2$, where $h_2 \in {\cal P}_2$, and
$h_1 + h_3 + h_2 \in {\cal P}_1$. Let us rewrite the image of $\DV_{12;43}$
in the form $h_1 + h_3 + h_2 = h_1 - (h_1 - (h_1 + h_3 + h_2))$. Since
$h_1 + h_3 \in {\cal P}_4$, and $h_1 -(h_1 + h_3) \in {\cal P}_3$, then
$h_1 -(h_1 + h_3)$ is the image of $h_1$ under the projection $\pi_4$
on ${\cal P}_3$. The subtraction of $h_2$ gives the
image of projection $\pi_2$
on ${\cal P}_1$. Hence $(h_1 - (h_1 + h_3 + h_2)) = \DV_{12;34}h_1$.
Consequently, $\DV_{12;43} = \Id - \DV_{12;34}$.
$\Box$

Lemmas 1-3 define generators of the representation of the group of
permutations for four subspaces in the group generated by identity and
$D = \DV_{12,34}$. Let us formulate corollaries from these lemmas
first for various composition mappings of the space
${\cal P}_1$ on itself.

$$
{\DV}_{12,43} = \Id - D; \, {\DV}_{13,24} = {\DV}_{14,23} =
(\Id - D^{-1})^{-1}, \,
{\DV}_{13,24} = (\Id - D)^{-1}, \, {\DV}_{14,32} = D^{-1}.
$$

It remains to consider the case when the images (and the preimages) of
composition mappings are nonisomorphic (say have different dimensions).
The corresponding operators in this case differs only by direct
summonds of bigger subspace which maps identically. Reduction to the
case of isomorphic subspaces is realized as follows. Let
$ \dim {\cal P}_2 > \dim {\cal P}_1$, and one wishes to compare
$\DV_{12;34}$ with $\DV_{21;43}$.
Under the mapping $\DV_{21;43}$ the subspace
$({\cal P}_2 \cap {\cal P}_4)$ maps on itself identically because its
vectors belong both to ${\cal P}_2$ and ${\cal P}_4$.
The subspace ${\cal P}_{13} = ({\cal P}_1 \oplus {\cal P}_3)$
is invariant also because both mappings
$\pi_1$ and $\pi_3$ are projections parallel to this subspace.
${\cal P}_{13}$ is invariant relative to $\DV_{12;34}$ also
because it contained both subspaces ${\cal P}_1$, and ${\cal P}_3$.
Hence, the question is reduced to the structure of restriction of mappings
$\DV_{12;34}$ and $\DV_{21;43}$ on ${\cal P}_{13}$.
In general case the intersection of subspaces ${\cal P}_2$ and
${\cal P}_4$ with ${\cal P}_{13}$ have the same dimension as
${\cal P}_1$. Hence, one can apply lemmas 1-3 in the space
${\cal P}_{13}$:
$$
{\DV}_{21,43} = \Id - {\DV}_{21,34} = \Id - {\DV}_{34,21} =
{\DV}_{34,12} = {\DV}_{12,34}.
$$

Let us show that operator cross-ratio defines cocycle with values in
operators.

\begin{lemma}
\label{lemma4}
Let two subspaces ${\cal P}_i, \; i=1,2$ be equivalent relative to the
group $U_{res}$ and three subspaces
${\cal Q}_j, \; j=1,2,3$ complete them to the whole space
(that is ${\cal P}_i \oplus {\cal Q}_j = H, \; i=1,2, \; j=1,2,3$).
Then
\beq
\DV({\cal P}_1{\cal Q}_1,
{\cal P}_2{\cal Q}_2) \DV({\cal P}_1{\cal Q}_3,{\cal P}_2{\cal Q}_1)
\DV({\cal P}_1{\cal Q}_2,{\cal P}_2{\cal Q}_3) = \mathop{\rm Id},
\label{2}
\eeq
\noindent i.e. the product of these three operators is the identity.
\end{lemma}

\Proof.  Consider the chain of mappings
$$
{\cal P}_1 \stackrel{{\cal Q}_2}{\longrightarrow}{\cal P}_2
\stackrel{{\cal Q}_1}{\longrightarrow}{\cal P}_1
\stackrel{{\cal Q}_1}{\longrightarrow}{\cal P}_2
\stackrel{{\cal Q}_3}{\longrightarrow}{\cal P}_1
\stackrel{{\cal Q}_3}{\longrightarrow}{\cal P}_2
\stackrel{{\cal Q}_2}{\longrightarrow}{\cal P}_1,
$$
\noindent which defines the left hand side of the formula~(\ref{2}).
The composition of the second and the third mappings, as well as the
composition of the fourth and the fifth ones, are identities. After
its reduction the remaining composition gets identity too.
$\Box$

Let us clarify the geometrical sense of the lemma 4.

Consider the canonical (tautalogical) bundle $\gamma$ on the manifold
$Gr(H)$, i.e. the bundle whose fiber at any point $W \in Gr(H)$ is the
linear space $W$. Introduce the following trivialization of $\gamma$.
Fix a point $W_+ \in SGr(H)$. The chart ${\cal U}_V$ on $\gamma$ is
defined by a plane $V \in H$ that complete $W_+$, i.¥. $W_+ \oplus V = H$,
and besides the projecting operator $\pi_V: H \to W_+$ parallel to the
space $V$ is assumed to be bounded, in other words, for the decomposition
$h = w + v$, where $h \in H, \, w \in W_+, \, v \in V$ the following
estimate with a constant $C$ has to be valid $||w|| \le C ||h||$ (the space
$V$ does not have "infinitesimally small" angles with $W_+$). The
coordinates of a point $(W, x) \in \gamma$, where $x \in W$, in the
chart ${\cal U}_V$ will be $(W, \pi_Vx) \in Gr(H) \times W_+$.
Let us calculate the transformation formulas of coordinates from a
chart ${\cal U}_{V_1}$ to that of ${\cal U}_{V_2}$.
Let coordinates of a point
$(W,x) \in \gamma$ in the chart ${\cal U}_{V_1}$ be
$(W, y)$, where $y \in W_+$. Then $x = \pi_{V_1}^{-1}y$, and the same
point has in the chart ${\cal U}_{V_2}$ coordinates
$(W, \pi_{V_2} \circ \pi_{V_1}^{-1} y)$. But $\pi_V = \pi_{V}^{-1}$,
since $\pi_V$ is a projector. Hence, the transition function is defined
by formula
$$
W_+ \stackrel{V_1}{\longrightarrow}W
\stackrel{V_2}{\longrightarrow}W_+,
$$
\noindent which coincides with the cross-ratio of four subspaces
${\DV}(W_+,V_2; W,V_1)$. By lemma 4, the transition from
${\cal U}_{V_1}$ to ${\cal U}_{V_2}$, further on, to
${\cal U}_{V_3}$, and finaly, back to ${\cal U}_{V_1}$ gives the cocycle
property~(\ref{2}).

So, the transition from a chart ${\cal U}_{V_1}$ to a chart
${\cal U}_{V_2}$ is defined by the transform ${\DV}(W_+,V_2;W,V_1)$
which acts on coordinates $x \in W$ as a linear fractional function
from operator coordinates of a plane $W$. These transforms are defined
on intersections of charts ${\cal U}_V$. But the set of these charts
does not cover all the Grassmann manifold. In contrast with the
finite-dimensional case, two isomorphic subspaces of
infinite-dimensional Hilbert manifolds does not have in general a
common complementary subspace $V$. The exact necessary and sufficient
conditions of existence a common complement for two given subspaces
was found in~(\cite{T}). To overcome this difficulty let us exchange
the cross ratio $\DV (W_+,V_2;W,V_1)$ by $\DV (W,V_1;W_+,V_2)$. In
view of lemma 1 the cross ratio remains the same but now it defines a
transformation of $W$ instead of $W_+$. So we can change $ W_+$ and
the corresponding charts cover all Grassmann manifold.

Hence, coordinate transformations of the canonical bundle $\gamma$
are endomorphisms $\DV (W,V_1;W_+,V_2)$ that can be regarded as
endomorphisms of $\gamma$ itself. In view of lemma 4 the cross ratio
defines on the Grassmannian a cocycle $\{\DV\}$ with coefficients in
the sheaf of endomorphisms of the canonical bundle $\gamma$, i.e.

$$
\{\DV\} \in H^1(SGr(H),\mathop{\rm End}(\gamma)).
$$

\smallskip

Following the scheme of Atijah \cite{A} and Turin \cite{Tu}
consider the principal bundle $P$ with the group $G$ corresponding to
the vector bundle $\gamma$.Take the exact sequence of bundles over
$SGr(H)$

\beq
0 {\longrightarrow} L \stackrel{f}{\longrightarrow} Q
\stackrel{g}{\longrightarrow} T {\longrightarrow} 0,
\label{3}
\eeq
where $T$ is the tangent bundle to $SGr(H)$; $Q$ is the bundle of
invariant tangent vector fields on $P$; $L$ is the bundle of Lie
algebras corresponded to left invariant vector fields on G
(vector fields tangent to fibers). The exact sequence (\ref{3})
defines an extension of $T$ by $L$. Classes of equivalent extensions
are in one-to-one correspondence with the elements of
$H^1(SGr(H),\mathop{\rm Hom}(T,L))$ --- the one-dimensional cohomology
group with coefficients in the sheaf $\mathop{\rm Hom}(T,L)$. As in the
finite-dimensional case, we say that the sequence (\ref{3}) is split
if there exists a homomorphism $h:T \to Q$ such that
$gh = \Id : T \to T$.

A connection in the principal bundle is a splitting of the
corresponding exact sequence. Let us denote by
$$
a(\gamma) \in  H^1(SGr(H),\mathop{\rm Hom}(T,L))
$$
an element corresponding to the extension (\ref{3}).
The coordinate transformations $\varphi_{ij} = u_i^{-1}u_j$ from a
chart $u_i: U_i \times G \to P|_{U_i}$ to a chart
$u_j: U_j \times G \to P|_{U_j}$ of the bundle $P$ are defined by the
cross ratio $\varphi_{ij} = \DV (W,V_i;W_+,V_j)$. The chart $u_i$
induces an isomorphism of tangent bundles and since it commute with
the action of $G$ we obtain the isomorphism
$\hat u_i: T_i \oplus L_i \to Q_i$. Here $T_i, Q_i, L_i$ are the
restrictions of the corresponding bundles to the neighbourhood $U_i$.
The sequence (\ref{3}) splits on this neighbourhood hence there exists
a lifting of the identity endomorphism of $T$ that gives an element
$a_i: T_i \to Q_i$, namely $a_i(t) = \hat u_i(t \oplus 0)$. Put
$a_{ij} = (a_j -a_i): T_{ij} \to Q_{ij}$. Then $\{a_{ij}\}$ represents
the cocycle $a(P)$.

Let $\Omega^1$ be the sheaf of germs of differential 1-forms on the
manifold $SGr(H)$. We have $\Omega^1 = \mathop{\rm Hom}(T,1)$. Hence
$H^1(SGr(H),\mathop{\rm Hom}(T,L)) = H^1(SGr(H),(L \times \Omega^1))$.
Thus $a()P \in H^1(SGr(H),(L \times \Omega^1))$. For compact K\"ahler
manifolds in view of isomorphism Dolbeault
$H^1(SGr(H),(L \times \Omega^1))$ correspondes to cohomologies of
$H^{1,1}$-type. It generates the ring of characteristic classes of the
corresponding bundle. In our case (the manifold $SGr(H)$ is noncompact)
we will by definition consider $a(P)$ as an analog of the generating
Chern class of the canonical bundle $\gamma$.

\bigskip
\bigskip

The group $U_{res}$ transitively acts on the big cell of Sato's
Grassmannian but not double transitively, as it is known even in
finite-dimensional case. It is natural to find invariants of pairs
of points relative to the group $U_{res}$. These invariants are
designed as in finite-dimensional case~\cite{Z} and are closely
related with the operator cross-ratio.

Let ${\cal A}; \, {\cal B} \in Gr_+(H)$. Its coordinates are
Hilbert-Shmidt operators $A$ and $B$ respectively.
The classe of operator
$\DV_{AB}:=(\Id+A^*A)^{-1}(\Id+A^*B)(\Id+B^*B)^{-1}(\Id+B^*A)$
will be called the operator angle between subspaces
${\cal A}$ and ${\cal B}$. Let us remark that the operator $(\Id+A^*A)$
is positive definite and invertible.

Two pairs of subspaces
${\cal S},\, {\cal T}$ and ${\cal P},\, {\cal Q}$ shall be called
equivalent relative to $U_{res}$, if there exists an element
$g \in U_{res}$, such that $g{\cal S} = {\cal P},\, g{\cal T} = Q.$

Two operators $V,W:H_+ \to H_-$ shall be called comparable if there
exist two unitary operators
$\alpha \in U(H_-)$ and
$\beta \in U(H_+)$, such that $W = \alpha V \beta.$

\begin{theorem}
\label{theorem2}
Two pairs of points ${\cal P},\, {\cal Q}$ and ${\cal S},\, {\cal T}$
of Sato's Grassmannian are equivalent iff the classes
$\DV_{PQ}$ and $\DV_{ST}$ coincide.
\end{theorem}

First we prove the following lemma.

\begin{lemma}
\label{lemma5}
The operator angle $\DV_{AB}$ between subspaces
${\cal A}$ and ${\cal B}$ is the operator of
the composite mapping
$$
{\cal A} \stackrel{{\cal B}^{\perp}}{\longrightarrow}{\cal B}
\stackrel{{\cal A}^{\perp}}{\longrightarrow}{\cal A}.
$$
\end{lemma}

\Proof. Let us decompose an element
$(h,Ah)$ of the subspace ${\cal A}$ relative to the basis
${\cal B},{\cal B}^{\perp}$. If $(x, Bx)$ is its projection on
${\cal B}$, then $(h,Ah) - (x, Bx)$ is orthogonal to ${\cal B}$, i.e.
$\langle (h-x, Ah-Bx), (y,By) \rangle = 0$ for any $y$ where
angle brackets stand for scalar product in $H$. It follows that
$x = (\Id+B^*B)^{-1}(\Id+B^*A)h$. Making the decomposition
relative to the basis
${\cal A},{\cal A}^{\perp}$, one obtains the element
$(z,Az)$ with
$z = (\Id+A^*A)^{-1}(\Id+A^*B)(\Id+B^*B)^{-1}(\Id+B^*A)h$.
It remains to apply the proposition 1.

$\Box$

{\bf Proof of the Theorem 2.}

Necessity follows from the lemma 5 in view of invariance of the
composite mapping

Sufficiency.

Since $U_{res}$ transitively acts on the big cell, let us transfer
${\cal S}$ into ${\cal P}$ and put the coordinate of the obtaining
subspace equal to zero. Let the coordinate of ${\cal T}$ became $W$
and the coordinate of ${\cal Q}$ became $V$. Then
$\DV_{PQ} = (\Id+V^*V)^{-1}, \, \DV_{ST} = (\Id+W^*W)^{-1}$.
Since these operators are similar then operators $V^*V$ and $W^*W$ are
similar also. These operators correspond to the quadratic forms
$\langle Vx,Vx \rangle$ and $\langle Wx,Wx \rangle$ correspondinly
on the space $H_+$.

To continue the proof we need the following lemma.

\begin{lemma}
\label{lemma6}

Let $V,W:H_+ \to H_-$ are Hilbert-Shmidt operators.
Two quadratic forms
$\langle Vx,Vx \rangle$ and
$\langle Wx,Wx \rangle$ on the space $H_+$ are similar iff operators
$V$ and $W$ are comparable.
\end{lemma}

Sufficiency.

Let $W = \alpha V \beta$; $\alpha ^* \alpha = \Id; \;
\beta ^* \beta = \Id.$ Then
$$
(W^*W) = \beta ^* V^* \alpha ^* \alpha V \beta =
\beta ^* (V^*V) \beta.
$$

Necessity.

Let $(W^*W) = \beta ^* (V^*V) \beta.$ Then $\ker W =
\beta (\ker V)$. Multiply $V$ from the right by the unitary operator
$\beta$ such that $\ker W = \ker V_1$, where $V_1 = V\beta$.
In view of proved sufficiency the operators
$V^*_1V_1$ and $W^*W$ are similar. Multiply $V_1$ from the left by
the unitary operator
$\gamma:H_- \to H_-$, which transfers the image of
$\mathop{\rm Im} V$ into $\mathop{\rm Im} W$.
Denote $\gamma V_1 = V_2$.
Operators $V_2^*V_2$ and $W^*W$ remain similar.

The operator $\tilde V:(\ker W)^{\perp} \to \mathop{\rm Im} W$, being
restriction of the operator
$V_2$, is invertible. Denote by $\tilde W$ the restriction of the
operator $W$ to the subspace $(\ker W)^{\perp}$.
The operators $\tilde V$ and $\tilde W$ are equivalent. Hence, there
exists a unitary transform $q$ of the subspace $(\ker W)^{\perp}$
such that

\beq
\tilde W^* \tilde W = q^* \tilde V^* \tilde V q.
\label{w}
\eeq
\noindent Take $\alpha_1 = q^*, \; \alpha_2 =
\tilde W q^* \tilde V^{-1}$. Then, firstly, $\alpha_2 \tilde V
\alpha_1 = \tilde W$, i.e. $\tilde V$ and $\tilde W$ are comparable,
and, secondly, $\alpha_2^* \alpha_2 = (\tilde V^*)^{-1}q \tilde W^*
\tilde W q^* \tilde V^{-1} = \Id$, in view of~(\ref{w}).
It follows that the operator
$\alpha_2$ is unitary. If we extend $\alpha_1$ and $\alpha_2$ by
identity to the kernel we obtain that $V$ and $W$ are comparable.
$\Box$

To complete the proof of the theorem it remains to note that the
stabilizer of the point
$0$ in $U_{res}$ are block diagonal operators
$$
\left ( \alpha_1 \quad 0 \atop 0 \quad \alpha_2 \right ),
$$
\noindent transferring coordinate $V$ of the subspace ${\cal Q}$
into the comparable operator. Hence, the pair of subspaces
${\cal P},{\cal Q}$ transforms into the pair ${\cal S},{\cal T}$.

$\Box$

\bigskip

{\bf Integrals of KP-Hierarchy}

\bigskip

While finding $\tau$-function~\cite{PS} one obtains an infinite matrix
which was upper-triangular with units on the main diagolal, except for a
finite-dimensional block. That was the reason for existing its
determinant. It is easily seen that it is the matrix of operator
cross-ratio. If we use the lemma 3 and change the order of four subspaces
we obtain the infinite matrix $A$ which will be "asymptotically
nilpotent" (i.e.  upper-triangular with zeros on the main diagolal)
except for a finite-dimensional block. Such matrix $A$ we shall call
almost asymptotically nilpotent. Almost asymptotically nilpotent
matrix and all its powers have a trace which gives us the possibility
to define $\zeta$-function of a cross-ratio.

The cross-ratio relates to four spaces
${\cal P}_1,{\cal P}_2,{\cal P}_3,{\cal P}_4 \in Gr_+$
the class of operators similar to
$\DV({\cal P}_1,{\cal P}_2,{\cal P}_3,{\cal P}_4)$. Each subspace
${\cal P}_i$ defines, due to Sato's correspondence, a solution of
KP-Hierarchy therewith the image of each ${\cal P}_i$ is obtained by
M\"obius transform. Hence, the cross-ratio of four spaces remains
invariant in the flow of KP. It means that it changes
isospectrally. Invariants
of almost asymptotically nilpotent matrix relative to isospectral
deformations are traces of its powers. It is natural to describe
these invariants by $\zeta$-function.

Let we have three solutions $W_i(t), \; i=1,2,3$
of Riccati equation on Sato's Grassmannian. Using Sato's
correspondence we obtain three solutions $L_i(t), \; i=1,2,3$
of KP-hierarchy on the space $E'$. Take any subspace $W$ and designe
the cross-ratio $\DV (W,W_1(t),W_2(t),W_3(t))$.
We obtain the operator, spectrum of which does not depend on $t$. Its
characteristic numbers gives the set of integrals and $\zeta$-function
remains invariant. It is convenient for calculations to take,
as some of $W_i(t)$, stationary solutions of the corresponding Riccati
equations.

\bigskip
\bigskip
\bigskip

{\bf Schwarz operator.}

\bigskip

Operator cross ratio allows us to define operator analog of Schwarz
derivative for curves taking values in Grassmann manifolds \cite{Z}.

Let $H$ be a Hilbert space equiped with a polarization
$H = H_+ \oplus H_-$. Consider the Grassmann manifold $Gr(H)$ defined
by the subspace $H_+$ and the dual Grassmannian $\tilde {Gr}(H)$
consisted of complementary subspaces. We will consider the case when
complementary subspaces are isomorphic to that of $Gr(H)$. In the
finite-dimensional case this means that dimensions of $H_+$ and $H_-$
coincide.

Consider a smooth family of subspaces of $H$ depending of
onedimensional parameter $s$ such that pairs
$(z(s),z(\sigma))$ for $s > 0$ and $\sigma \le 0$ give a polarization
of the space $H$. Our nearest aim is to build an analog of the
Schwarz derivative of the curve $z(s)$ at $s=0$.

Take $s_2 < 0 < s_1 < s_3$. Then pairs of spaces
$(z(s_2),z(s_1))$ and $(z(0),z(s_3))$ give polarizations of $H$.
From now on we will for brevity write $z_i, z'_i ...$ instead of
$z(s_i), z'(s_i) ...$; values of $z(0), z'(0) ...$ will be denoted by
$z, z' ...$ (without indices). Consider the mapping $f$ depending on
four parameters and giving by the cross ratio

$$
f(s_2,s_1,0,s_3) = \DV (z_2,z_1;z,z_3) =
(z_2 - z_1)^{-1}(z_1 - z)(z - z_3)^{-1}(z_3 - z_2).
$$

Let $s_2 \to 0$. At $s_2=0$ the cross ratio $f$ gives the identity
mapping. The derivative of $f$ at $s=0$ is equal to

$$
\frac {\partial f}{\partial s_2}(0,s_1;0,s_3) =
- (z-z_1)^{-1}z' + (z-z_3)^{-1}z'.
$$

Now let $s_3 \to s_1$. At $s_3 = s_1$ the function $f$ and its
derivative relative to $s_2$ are equal to zero. We have

\beq
\frac {\partial^2 f}{\partial s_3\partial s_2}(0,s_1;0,s_1) =
- (z-z_1)^{-1}z'_1(z-z_1)^{-1}z'.
\label{z}
\eeq

The right hand side of (\ref{z}) is defined at $s_1 > 0$ and has a
singularity at $s_1=0$. Consider the situation when subspaces of a
polarization tend to be bound together. Find the asymptotic of (\ref{z})
when $s_1 \to 0.$

$$
\begin{array}{l}
z'_1 = z' + s_1z'' + \frac {s_1^2z'''}{2} + o(s_1^2), \\
z_1 - z = s_1z'\left ( \Id + \frac {s_1(z')^{-1}z''}{2} +
\frac {s_1^2(z')^{-1}z''}{6} + o(s_1^2) \right ), \\
(z_1 - z)^{-1} = s_1^{-1}\left ( \Id - \frac {s_1(z')^{-1}z''}{2} -
\frac {s_1^2(z')^{-1}z'''}{6} +
\frac {s_1^2((z')^{-1}z'')^2}{4} + o(s_1^2) \right ).
\end{array}
$$

By substituting these expressions for (\ref{z}) we obtain
$$
\frac {\partial^2 f}{\partial s_3\partial s_2}(0,0;0,0) =
(s_1)^{-2}\left ( \Id + \frac {s_1^2(z')^{-1}z'''}{6} -
\frac {s_1^2((z')^{-1}z'')^2}{4} + o(s_1^2) \right ).
$$

In accordence with this formula define the differential operator

\beq
S(z) = (z')^{-1}z''' - \frac {3}{2}((z')^{-1}z'')^2.
\label{s}
\eeq

This expression is an analog of the Schwarz derivative. In the
finite-dimensional case it was introduced and explored in \cite{Z}.
The above given deduction shows that the expression (\ref{s}), as a
cross ratio with the help of which it is built, defines the same class
of conjugate operators independently of M\"obius transforms in the
ambient space. Let us show that (\ref{s}) is closely relates with
Hamiltonian systems and Riccati equations.

Let
\beq
\left \{
\begin{array}{l}
q' = A(t)q + p  \\
p' = - B(t)q - A^{\ast}(t)p
\end{array}
\right .
\label{4}
\eeq
be a linear Hamiltonian system with the Hamiltonian

\beq
H = \frac {1}{2}(<p,p> + 2<p,Aq> + <Bq,q>),
\label{5}
\eeq
where $A(t)$ and $B(t)$ are $n \times n$-matrices; the matrix $B(t)$
is symmetric; angle brackets stand for scalar product in $\R^n$.
Hamiltonian (\ref{5}) corresponds to minimization problem of the
quadratic functional

\beq
\int_{t_0}^{t_1}(<q',q'> - 2<Aq,q'> + <(A^{\ast}A-B)q,q>)dt.
\label{6}
\eeq
The identity coefficient of the first summond is the result of
reduction of variational problems with the strong Legendre condition.
The system (\ref{4}) is equivalent to the Euler equation for the
functional (\ref{6}):

\beq
q'' + (A^{\ast} - A)q' + (B - A' - A^{\ast}A)q = 0.
\label{7}
\eeq

We will consider fundamental systems of solutions to (\ref{4})
and hence $p$ and $q$ are $n \times n$-matrices also.

The coeffitions matrix of the system (\ref{4}) is symplectic thus
(\ref{4}) defines a flow on the Lagrange-Grassmann manifold $\Lambda$
of the space $(p,q)$. Local coordinates of points of $\Lambda$ are
given by matrices $W = pq^{-1}$. The evolution of coordinates $W$ is
described by Riccati equation
\beq
W' = (-Bq-A^{\ast}p)q^{-1} - pq^{-1}(Aq+p)q^{-1} = -B-A^{\ast}W-WA-W^2.
\label{8}
\eeq

\smallskip

Let us return to our analog of Schwarz operator.

\beq
\begin{array}{l}
S(z(t)) = [(z'(t))^{-1}z'']' -
\frac {1}{2} [(z'(t))^{-1}z'']^2 = \\
(z'(t))^{-1}z''' -
\frac {3}{2}[(z'(t))^{-1}z'']^2.
\end{array}
\label{9}
\eeq

It is convenient to consider $t$ as changing on a projective line
(real or complex) or on a Rimanian surface on which acts M\"obius
group of linear fractional mappings.

It was shown in \cite{Z} that the equivalent class of the image of
Schwarz derivative $S(z(\cdot))$ is invariant relative to the
M\"obius group. Namely if one denote by $M$ a M\"obius transform
$M : z \mapsto (C_1z + C_2)(C_3z + C_4)^{-1}$, where
$C_i$ is $(n \times n)$-matrices, then there exists a matrix $K(t)$
such that $S(M(z(t))) = K(t)S(z(t))K^{-1}(t)$. In other words,
M\"obius transforms of preimage of the Schwarz operator lead to
isospectral change of the image.

Let us describe a connection of the operator $S$ with the Hamiltonian
system (\ref{4}) or (it is the same) with the Riccati equation
(\ref{8}). Suppose that the matrix $A$ is symmetric. We find a
connection between a solution to Riccati equation $W$ and the function
$z$ given by the formula

\beq
W = - \frac {1}{2} [(z'(t))^{-1}z''] - A.
\label{10}
\eeq

We have
$$
\begin{array}{l}
W' = - \frac {1}{2} (z'(t))^{-1}z''' +
\frac {1}{2} (z'(t))^{-1}z''(z'(t))^{-1}z'' - A'; \\

W^2 = \frac {1}{4} (z'(t))^{-1}z''(z'(t))^{-1}z'' +
\frac {1}{2} (z'(t))^{-1}z''A +
\frac {1}{2} A (z'(t))^{-1}z'' + A^2
\end{array}
$$
Hence
$$
\begin{array}{l}
W' + W^2 = - \frac {1}{2} (z'(t))^{-1}z''' +
\frac {3}{4} [(z'(t))^{-1}z'']^2(z'(t))^{-1}z'' - A' \\

\frac {1}{2} (z'(t))^{-1}z''A +
\frac {1}{2} A (z'(t))^{-1}z'' + A^2
\end{array}
$$
In view of (\ref{10}) the last formula gives
$$
\begin{array}{l}
W' + W^2 = - \frac {1}{2} S(z) - A' + A^2 - \\
\frac {1}{2} (2W+A)A -
\frac {1}{2} A (2W+A)
\end{array}
$$
or
\beq
W' + W^2 + WA + AW = - \frac {1}{2} S(z) - A'.
\label{11}
\eeq

Let us call
\beq
S(z(t)) = [(z'(t))^{-1}z'']' -
\frac {1}{2} [(z'(t))^{-1}z'']^2 = 2B(t) - A(t)
\label{12}
\eeq
{\sl Schwarz equation} for Hamiltonian system (\ref{4}) or (which is
the same) for Riccati equation (\ref{8})).

We proved the following.

\begin{theorem}

If $z(\cdot)$ is a solution to Schwarz equation (\ref{12}), then
$W(\cdot)$, defined by (\ref{10}), is a solution to Riccati equation
(\ref{8}).

If $W(\cdot)$ is a solution to Riccati equation (\ref{8}), then any
function $z(\cdot)$, being a solution of the linear
relative to $z(\cdot)$ equation (\ref{10}),
is a solution to Riccati equation (\ref{12}).
\end{theorem}

\bigskip
\bigskip

I am deeply gratefull to A.N.Parshin and A.A.Agrachev for stimulating
discussions and to Scuola Internazionale Superiore di Studi Avanzati
for the hospitality and excellent conditions for work.

\end{document}